\documentclass{article}
\usepackage{amsmath}
\usepackage{amsthm}
\usepackage{graphicx}

\newtheorem{proposition}{Proposition}[section]

\newtheorem{corollary}{Corollary}
\newtheorem{Conjecture}{Conjecture}
\newtheorem{Lemma}{Lemma}

\begin{document}
\title{Ehrenfest Wind-Tree Model is Dynamically Richer Than the Lorentz Gas}
\author{H. Attarchi, M. Bolding, and L. A. Bunimovich}
\maketitle
\begin{abstract}
	We consider a physical Ehrenfests' Wind-Tree model where a moving particle is a hard ball rather than (mathematical) point particle. We demonstrate that a physical periodic Wind-Tree model is dynamically richer than a physical or mathematical periodic Lorentz gas. Namely, the physical Wind-Tree model may have diffusive behavior as the Lorentz gas does, but it has more superdiffusive regimes than the Lorentz gas. The new superdiffusive regime where the diffusion coefficient $D(t)\sim(\ln t)^2$ of dynamics seems to be never observed before in any model.
\end{abstract}
\section{Introduction}
Consider the uniform motion of a point particle on the Euclidean plane in an array of rectangular or rhombus scatterers. When the particle hits the boundary of a scatterer it gets elastically reflected. This model was introduced in the celebrated paper by Paul and Tatyana Ehrenfest \cite{Ehrenfest} as a simple mechanical (dynamical) model for diffusion. In their model scatterers were rhombuses and this system itself was called a Wind-Tree model, where ``wind" stays for the particle and ``tree" for a scatterer \cite{Ehrenfest}. The Wind-Tree model was extensively studied by physicists \cite{BiaRon09,DetCohBei99,Gal69,H-C69,VanHau72,WooLad71}.\par
It turned out that this system is not a good model for diffusion. Indeed the Wind-Tree model is a billiard dynamical system and, because the boundary of the corresponding billiard table consists of straight segments, its dynamics is similar to billiards in polygons. It is well known that billiards in polygons have zero Kolmogorov-Sinai entropy and cannot generate dynamical chaos which is a necessary condition for demonstrating stochastic behavior and, in particular, diffusion. The Wind-Tree model was actively studied in statistical mechanics after appearance of powerful computers and found lacking of diffusion \cite{De13,DeHuLe14,H-C69,VanHau72,WooLad71}. Instead, the role of the basic simplest mechanical model of diffusion was very successfully played by the celebrated Lorentz gas where the scatterers are circles rather than rhombuses \cite{Bleher92,Bun81,BSC91}. Therefore these billiard systems belong to the class of the most chaotic billiards. In fact, Lorentz gas is unfolding of a Sinai billiard \cite{Si70}. Notably, the Lorentz gas was introduced as a model of electronic gas in metals which occurred to be completely irrelevant. Likewise, the Ehrenfests' Wind-Tree model is considered to be irrelevant as a model for diffusion and moved from statistical mechanics to pure mathematics where it is very popular now \cite{AvHu,De13,DeHuLe14,FrUl14,HarWeb80,HoHuWe13,HLT11,HuWe13,RaTr12}.\par
Our goal in this paper is to demonstrate that in fact, the Wind-Tree model is a good model to study diffusion and its dynamics is even richer than every bodies favorite Lorentz gas. Although we study only the physical periodic Wind-Tree model, it is quite clear that the non-periodic Wind-Tree model is probably dynamically richer than the (physical and mathematical) Lorentz gas with the same configuration in the plane of the centers of scatterers, but at least it has all the regimes of diffusion which the non-periodic Lorentz gas does.\par
A key observation is that a physical Wind-Tree model becomes a semi-dispersing billiard for any positive radius $r$ of a hard ball (disk in $\mathbf{R}^2$), whereas a physical Sinai billiard always remains Sinai billiard for any $r>0$. Therefore, in case of a bounded free path (finite horizon) both periodic Lorentz gas and periodic Wind-Tree model demonstrate diffusive behavior. However, if a free path is unbounded then dynamics of both periodic systems becomes superdiffusive and it is where Wind-Tree model overpasses Lorentz gas. If the configuration of scatterers is periodic then the particle may have unbounded free path only in strips on plane bounded by two parallel lines. Such strips traditionally are called corridors. In periodic Lorentz gas with unbounded free path, there is only one type of corridors, while in the periodic Wind-Tree model there are two types of corridors. Presence of a corridor of the first type results in the $(\ln t)$ growth of the diffusion coefficient while corridors of the second type make it grow as $(\ln t)^2$.\par
The structure of the paper is the following. In Section 2 we introduce the necessary notations and study some properties of the Ehrenfests' Wind-Tree model with infinite horizon. Section 3 and 4 deal with the calculation of the tail of the distribution of the free motion vectors (displacement) in two different types of corridors. In Section 5 an estimation of their correlations is given and the limit distributions of properly normalized free motion vectors in both discrete and continuous-time dynamics are discussed. 
\section{Ehrenfests' Wind-Tree Model}
\subsection{Configuration Space}\label{2.1}
\begin{figure}[h]
	\centering
	\includegraphics[width=10cm]{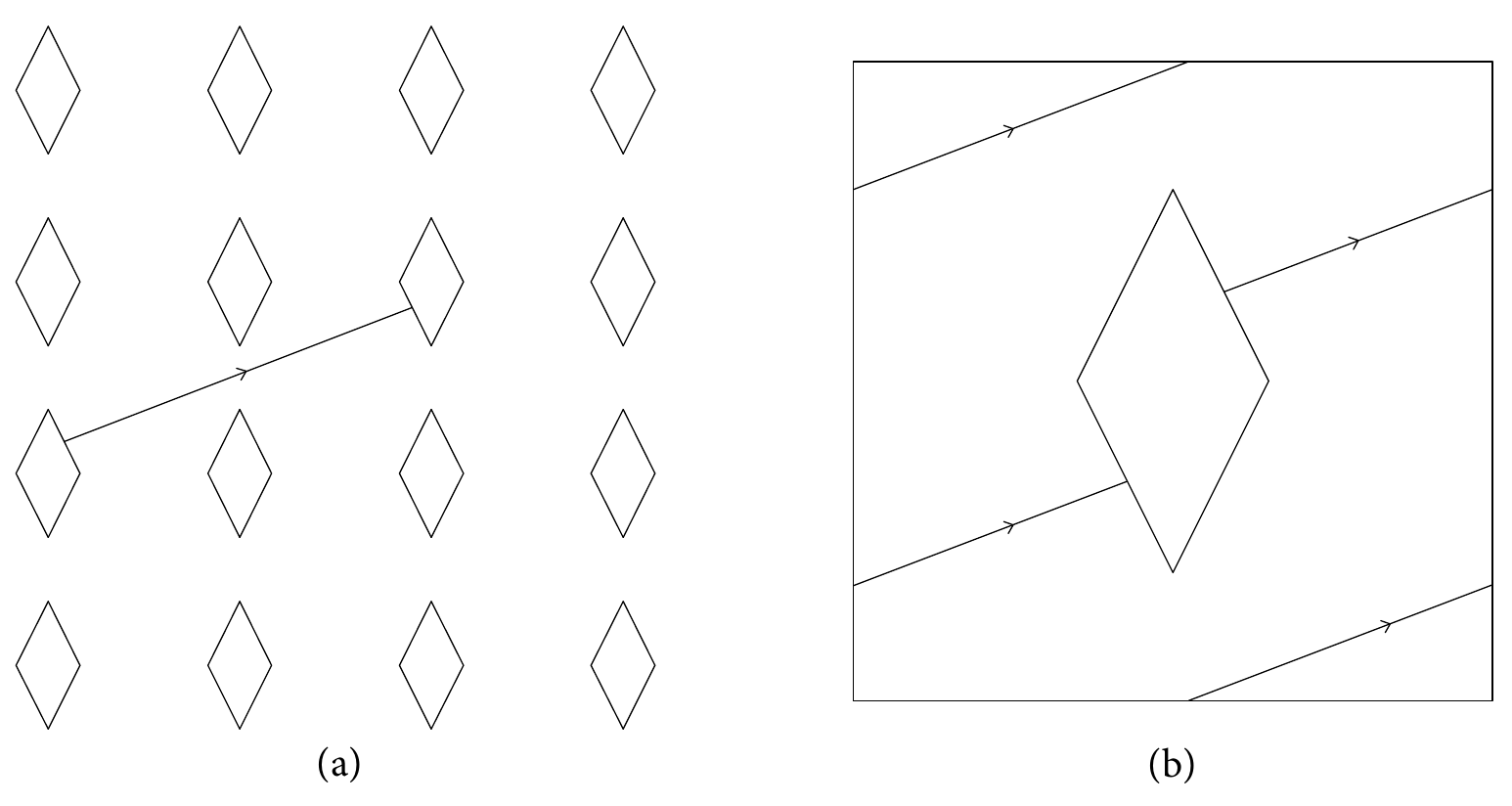}\caption{(a) Periodic Ehrenfests' Wind-Tree model. (b) Ehrenfests' Wind-Tree model in torus.}\label{fig0}
\end{figure}

The basic scatterers $S$ in Ehrenfests' Wind-Tree model are rhombuses. Each scatterer is determined by two parameters $\theta$ and $a$, where $2\theta$ is the acute angle of this scatterer and $a$ is its side length. Moreover, we assume that the diagonals of scatterers are parallel to $x$ and $y$-axes in $\mathbf{R}^2$ and acute angles are top and bottom angles of them. Then, in the periodic Ehrenfests' Wind-Tree model there is a periodic configuration of these basic scatterers with centers at the points with integer coordinates in $\mathbf{R}^2$ (Fig. \ref{fig0}a). Thus, a periodic Wind-Tree model is an unfolding of a billiard in a torus $\mathbf{T}^2$ with the scatterer $S$ where both are centered at the origin (Fig. \ref{fig0}b).\par
A configuration of scatterers is said to have a finite horizon if there exists $L>0$ such that any straight segment of the length bigger than $L$ in $\mathbf{R}^2$ intersects at least one scatterer (see Fig. \ref{fig1}b where we have two types of rhombuses). This means that the length of the free motion of the particle is bounded. Otherwise, the configuration of scatterers has an infinite horizon. When we have a configuration of scatterers with an infinite horizon, an open strip of parallel straight lines which intersects no scatterers is called a corridor. Moreover, a corridor has two opposite directions defined by the directions of parallel lines that bound the corridor. We will call two corridors equivalent if they have the same directions.\par
\subsection{Corridors in Ehrenfests' Wind-Tree Model}
\begin{figure}[h]
	\centering
	\includegraphics[width=12cm]{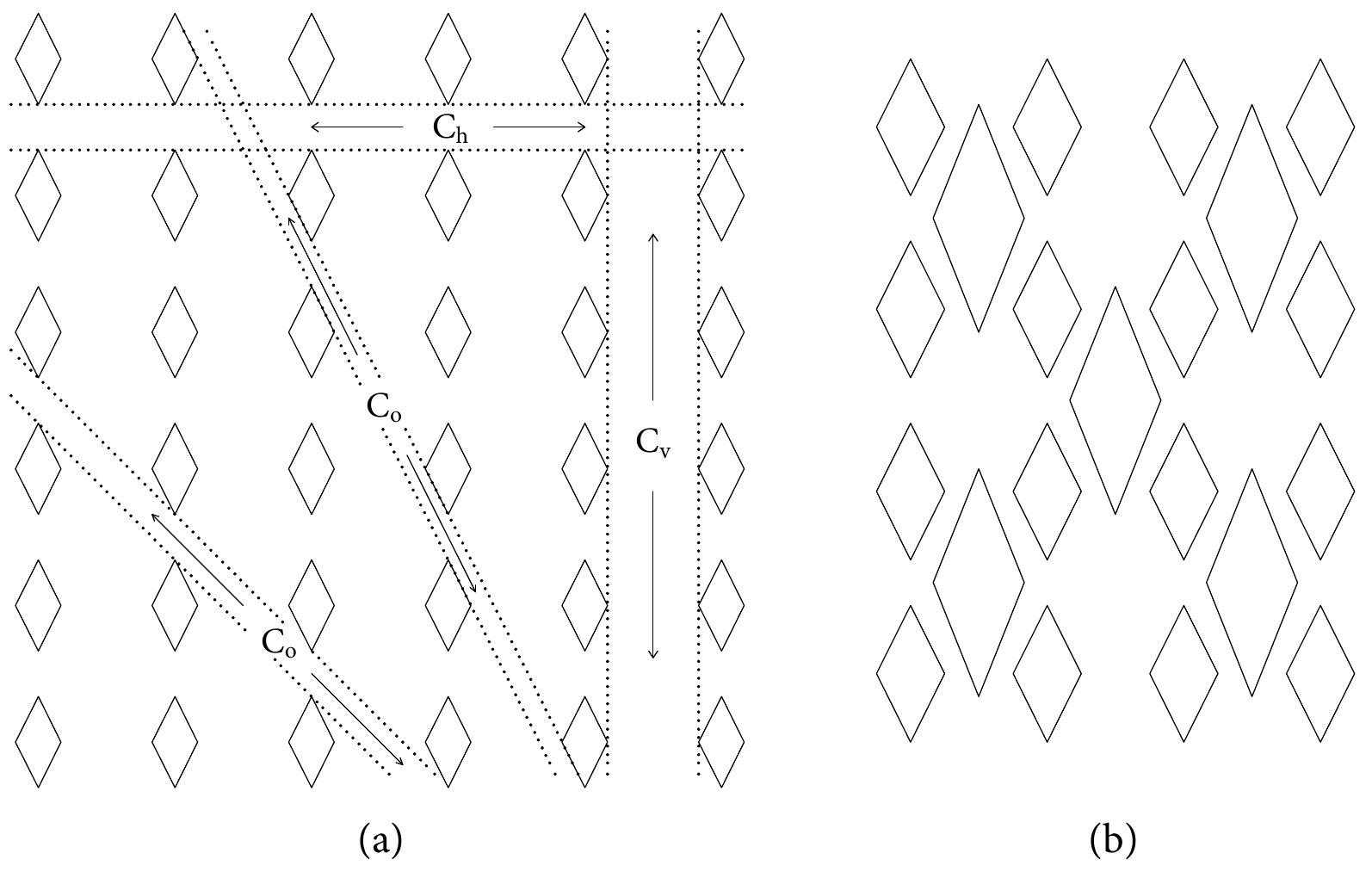}\caption{(a) Periodic Ehrenfests' Wind-Tree model with corridors: $C_h$, $C_v$, and $C_o$. (b) A periodic Ehrenfests' Wind-Tree model with a finite horizon and two types of rhombuses.}\label{fig1}
\end{figure}

The width of a corridor is defined by the distance between its boundary lines. The periodic Ehrenfests' Wind-Tree model, described in section \ref{2.1}, can have two types of corridors:
\begin{itemize}
	\item Type I: Corridors with boundaries touching some vertices of scatterers.
	\item Type II: Corridors with boundaries containing some edges of scatterers.
\end{itemize}
When $2a\cos\theta<1$, there exist horizontal and vertical corridors denoted by $C_h$ and $C_v$, respectively, and we call them corridors of type I (Fig. \ref{fig1}a). Their widths $d_h$ and $d_v$ are equal to:
$$d_h=1-2a\cos\theta,\hspace{1cm}d_v=1-2a\sin\theta,$$
respectively. In addition, for some choices of parameters $a$ and $\theta$ there exist oblique corridors which could be either of type I or II. We will denote the oblique corridors by $C_o$ and their width by $d_o$. (Fig. \ref{fig1}a)\par
\begin{figure}[h]
	\centering
	\includegraphics[width=5cm]{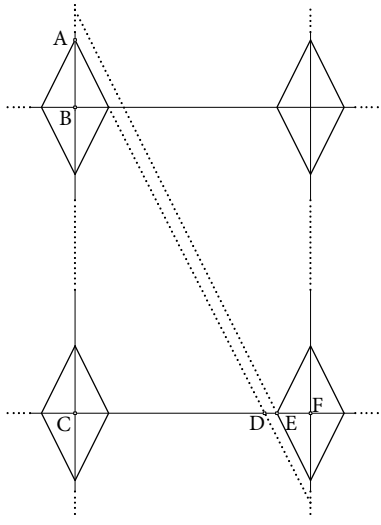}\caption{A passage through a corridor of type II between two adjunct columns of rhombuses.}\label{figlemma}
\end{figure}

\begin{Lemma}\label{lem1}
	There is an oblique corridor of type II if and only if $\tan\theta=\frac{m}{n}$ and 
	\begin{equation}\label{a<}
		a<\frac{\cos\theta+\sin\theta-\lceil\frac{n}{m}\rceil\sin\theta}{\sin(2\theta)}=\frac{(m+n-\lceil\frac{n}{m}\rceil m)\sqrt{m^2+n^2}}{2mn},
	\end{equation}
	where $0<m\leq n$ are integers. Moreover,
	$$d_o=\sin\theta+\cos\theta-\lceil\frac{n}{m}\rceil\sin\theta-a\sin(2\theta)=\frac{m+n-\lceil\frac{n}{m}\rceil m}{\sqrt{m^2+n^2}}-\frac{2mna}{m^2+n^2}.$$
\end{Lemma}
\proof
From the definition of corridors of type II and $\theta$, it is obvious that $\tan\theta\leq1$ is a rational number. To calculate the width of these corridors, we need to consider a passage by the particle through a corridor between two adjunct columns of rhombuses. According to Figure \ref{figlemma}, 
$$|CF|=1,\hspace{0.5cm}|EF|=a\sin\theta,\hspace{0.5cm}\tan\theta=\frac{|CD|}{|AB|+|BC|}=\frac{|CD|}{a\cos\theta+(\lceil\frac{n}{m}\rceil-1)}.$$
\par
\noindent Therefore,
$$|DE|=1-|EF|-|CD|=1-a\sin\theta-\left(\lceil\frac{n}{m}\rceil-1+a\cos\theta\right)\tan\theta.$$
Finally,
$$d_o=|DE|\cos\theta=\sin\theta+\cos\theta-\lceil\frac{n}{m}\rceil\sin\theta-a\sin(2\theta)$$
$$=\frac{m+n-\lceil\frac{n}{m}\rceil m}{\sqrt{m^2+n^2}}-\frac{2mna}{m^2+n^2}.$$
Solving the inequality $d_o>0$ for parameter $a$, we obtain:
$$a<\frac{\cos\theta+\sin\theta-\lceil\frac{n}{m}\rceil\sin\theta}{\sin(2\theta)}=\frac{(m+n-\lceil\frac{n}{m}\rceil m)\sqrt{m^2+n^2}}{2mn}.$$
These calculations also show that if  $\tan\theta=\frac{m}{n}$ and the inequality (\ref{a<}) are satisfied then $d_o>0$, where $d_o$ is the width of the oblique corridor of type II.\qed\par
Let $\alpha$ be the angle between the axis of an oblique corridor of type I and the positive $y$-axis such that $-\frac{\pi}{2}<\alpha<\frac{\pi}{2}$. It is easy to see that $\alpha\notin\{0,\pm\theta\}$.
\begin{Lemma}\label{lem2}
	The following relations hold for an oblique corridor of type I.\par
	\begin{enumerate}
		\item	$|\tan(\alpha)|=\frac{m}{n}$ for some integers $n>0$ and $m>0$.
		\item If $|\tan(\alpha)|<\tan\theta$, then
		$$d_o=\cos\alpha+\sin\alpha-\lceil\frac{n}{m}\rceil\sin\alpha-2a\sin\theta\cos\alpha=\frac{\left(n+m-m\lceil\frac{n}{m}\rceil-2an\sin\theta\right)}{\sqrt{n^2+m^2}},$$
		where $a<\frac{n+m-m\lceil\frac{n}{m}\rceil}{2n\sin\theta}$.
		\item If $\tan\theta<|\tan(\alpha)|\leq1$, then
		$$d_o=\cos\alpha+\sin\alpha-\lceil\frac{n}{m}\rceil\sin\alpha-2a\sin\alpha\cos\theta=\frac{\left(n+m-m\lceil\frac{n}{m}\rceil-2am\cos\theta\right)}{\sqrt{n^2+m^2}},$$
		where $a<\frac{n+m-m\lceil\frac{n}{m}\rceil}{2m\cos\theta}$.
		\item If $1<|\tan(\alpha)|$, then
		$$d_o=\cos\alpha+\sin\alpha-\lceil\frac{m}{n}\rceil\cos\alpha-2a\cos\theta\sin\alpha=\frac{\left(n+m-n\lceil\frac{m}{n}\rceil-2am\cos\theta\right)}{\sqrt{n^2+m^2}},$$
		where $a<\frac{n+m-n\lceil\frac{m}{n}\rceil}{2m\cos\theta}$.
	\end{enumerate}
\end{Lemma}
\proof
The proof is omitted since it is very similar to the proof of Lemma \ref{lem1}.\qed
\subsection{Physical Ehrenfests' Wind-Tree Model}\label{physical}
Physical billiards (i.e. billiards with a moving hard ball like a ball in real billiard as well as in all real systems modeled by billiards) were introduced and studied in \cite{Bun19}. A ball there and in the present paper is assumed to be a smooth hard ball. In case of a rough ball, it acquires rotation after collision at any point of the boundary \cite{CFZ18}.\par
It is easy to see that dynamics of a physical Lorentz gas is the same for any radius $r$ of the moving physical particle (of course unless $r$ becomes so big that the particle gets stuck in some subset of the plane). Indeed the scatterers remain to be circles but their radius increases by $r$. A totally different situation occurs in the Wind-Tree model where boundaries of scatterers acquire dispersing components (arcs of a circle of radius $r$). To represent the dynamics of a (hard) homogeneous spherical particle of radius $r>0$, it is enough to follow the motion of its center. It is easy to see that the center of particle moves in the smaller billiard table, which one gets by moving any point $q$ of the boundary by $r$ to the interior of the billiard table along the internal normal vector $n(q)$.\par
\begin{figure}[h]
	\centering
	\includegraphics[width=10cm]{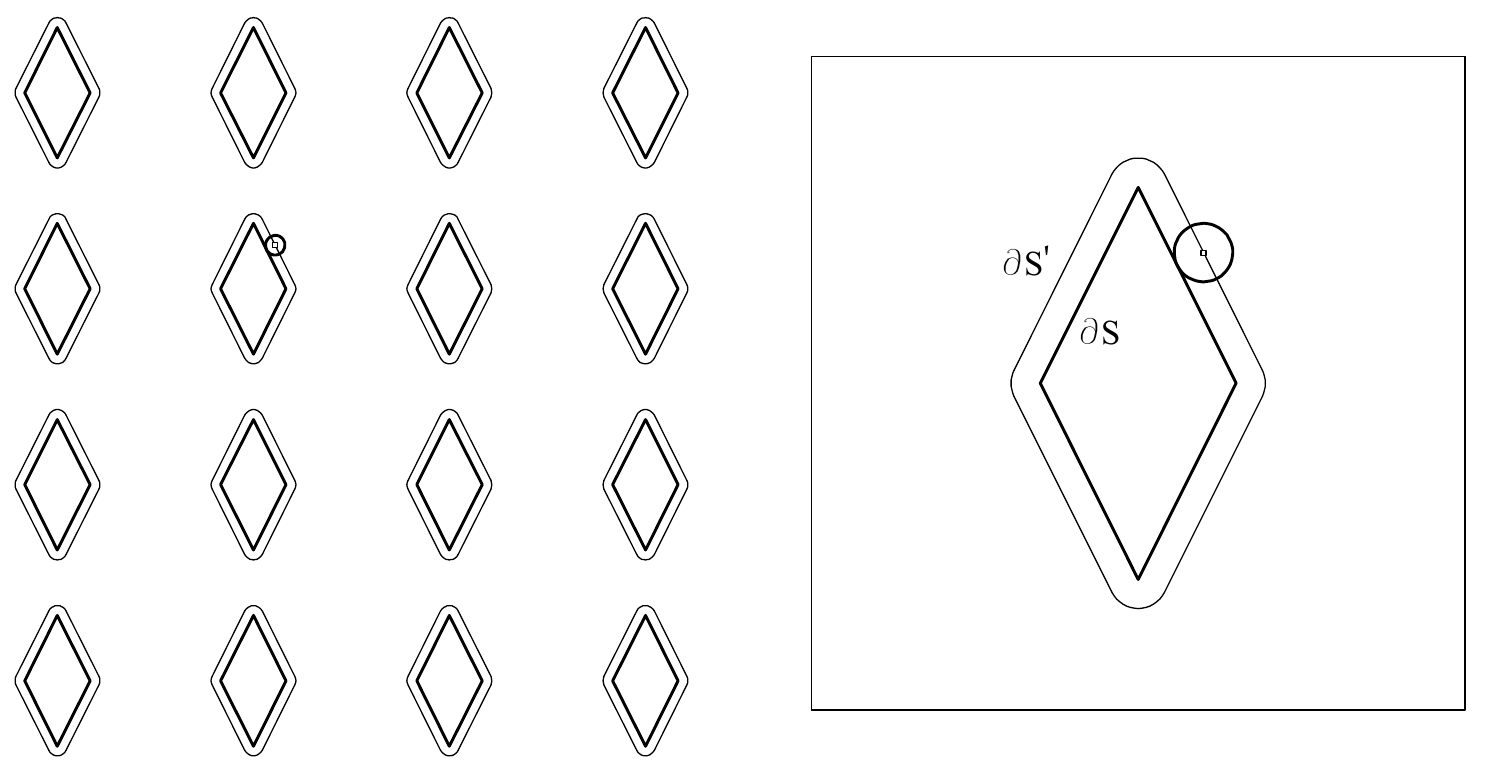}\caption{Physical Ehrenfests' Wind-Tree model and its equivalent mathematical billiard.}\label{fig11}
\end{figure}

In what follows, we will consider the mathematical billiard equivalent to the motion of a physical particle (disk) of radius $r>0$ in the periodic Ehrenfests' Wind-Tree model. We will use notations $S'$ and $\partial S'$ for the scatterer and its boundary in the equivalent mathematical billiard, respectively. (See Fig. \ref{fig11})\par
There are still two types of corridors in the equivalent mathematical billiard of the physical periodic Wind-Tree model:
\begin{itemize}
	\item Type I: Corridors with boundaries tangent to the dispersing components of scatterers.
	\item Type II: Corridors with boundaries containing the flat (neutral) components of scatterers.
\end{itemize}
Corridors remain if $r$ is small enough. More precisely, if $C$ is a corridor of type I (or II) in a physical periodic Ehrenfests' Wind-Tree model with the width $d$ and the physical particle has radius $r>0$ such that $d>2r$, then the equivalent mathematical billiard will have a corridor of type I (or II) in the same directions of $C$ with the width $d-2r$.\par
In the sequel, we will use the same notations $d_h$, $d_v$, and $d_o$ to denote the width of horizontal, vertical, and oblique corridors, respectively, in the equivalent mathematical billiard.
\subsection{Dynamics of the Physical Ehrenfests' Wind-Tree model}\label{dist1}
In this section, we introduce phase space of the system in the folded configuration space in torus.\par
Denote by $x(t)$ the position of the center of the physical particle at time $t$. We will characterize reflections by two quantities: $s$ the coordinate of the reflection point on $\partial S'$ with respect to the natural arc length $|.|$ in $\mathbf{R}^2$, and $\varphi$ the angle of reflection (the angle with the sign between the velocity vector after the reflection and the outer normal vector at the reflection point). Therefore, the phase space is:
$$\Lambda=\{X=(s,\varphi)\ |\ 0\leq s<|\partial S'|,\ \ -\frac{\pi}{2}\leq\varphi\leq\frac{\pi}{2}\}.$$
Now, let $T:\Lambda\rightarrow\Lambda$ be the Poincare section billiard map of the system such that $T(X_n)=T(s_n,\varphi_n)=(s_{n+1},\varphi_{n+1})=X_{n+1}$. This map preserves the Liouville measure:
\begin{equation}\label{mu}
	\mu(d\varphi ds)=Z^{-1}\cos(\varphi)d\varphi ds,
\end{equation}
where,
$$Z=\int_{\partial S'}\int_{-\pi/2}^{\pi/2}\cos(\varphi)d\varphi ds=2|\partial S'|,$$
on $\Lambda$. Moreover, the expected value of the function $F(X)$ is defined by:
$$\langle F(X)\rangle:=\int_\Lambda F(X)\mu(d\varphi ds)=Z^{-1}\int_{\partial S'}\int_{-\pi/2}^{\pi/2}F(X)\cos(\varphi)d\varphi ds,$$
where $X=(s,\varphi)$.\par
Let $x_n=x(X_n)\in\mathbf{R}^2$ be the position of the center of the physical particle in the discrete system at $n^{th}$ reflection. We will denote the segment (link) of the trajectory after $n^{th}$ reflection by $[x_n,x_{n+1}]$ and the corresponding vector of this free motion (displacement) by:
$$r(X_n):=x(TX_n)-x(X_n)=x(X_{n+1})-x(X_n)=x_{n+1}-x_n.$$
This implies,
\begin{equation}\label{x_n}
	x_n-x_0=x(T^nX_0)-x(X_0)=\sum_{i=0}^{n-1}r(T^iX_0).
\end{equation}
According to (\ref{x_n}), the problem of studying the statistical properties of the vector $x_n-x_0=x(T^nX_0)-x(X_0)$, when $n\rightarrow\infty$, is reduced to the problem about statistical properties of the free motion vector $r(X)$ with respect to the billiard map $T$. Let $\nu$ be the probability distribution of the free motion vector $r(X)$ with respect to Liouville measure $\mu$.
\begin{corollary}\label{nu-sym}
	The distribution $\nu$ is symmetric.
\end{corollary}
\proof
It is analogous to the proof of Proposition 4.1 of \cite{Bleher92}.\qed\par
In the following sections, we will study some properties of the distribution $\nu$ and of the second moment $\langle(r(X),r(T^nX))\rangle$, where $(.,.)$ denotes the standard inner product in $\mathbf{R}^2$ .
\section{Asymptotics of $\nu$ in Corridors of Type I}\label{sec3}
\begin{figure}[h]
	\centering
	\includegraphics[width=10cm]{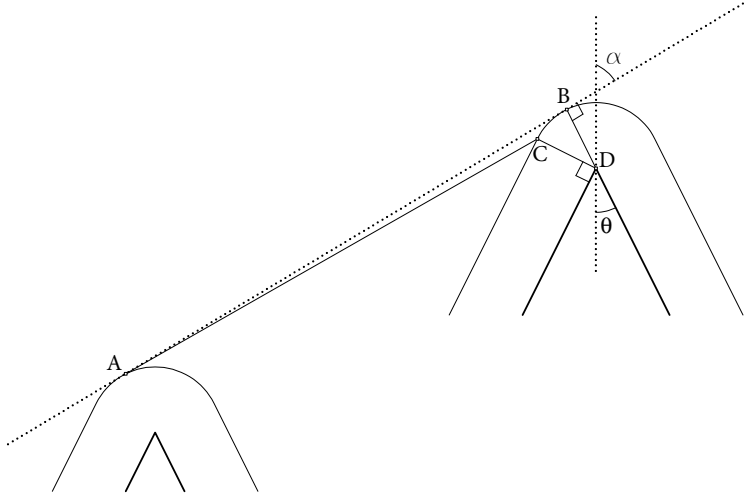}\caption{A Segment with an endpoint on a dispersing part tangent to a corridor of type I.}\label{fig3+}
\end{figure}

From Section \ref{physical}, we know that the boundary of a corridor of type I is tangent to the dispersing components of $\partial S'$. In the following lemma we will show that for a long enough segment $[x_0,x_1]$ in a corridor of type I, the endpoints $x_0$ and $x_1$ are on those dispersing parts of $\partial S'$ which are tangent to the boundary of the corridor.
\begin{Lemma}\label{lem3}
	Consider a finite segment $[x_0,x_1]$ of length $|r(X_0)|=|x_1-x_0|=L$ in a corridor of type I with the width $d$. There exists a constant $L_0>0$ such that if $L>L_0$, then the points $x_0$ and $x_1$ are on dispersing parts of $\partial S'$ which are tangent to the boundary of that corridor on its opposite sides.
\end{Lemma}
\proof
Let $\alpha$ be the angle between the direction of the corridor of type I and the positive $y$-axis as it is shown in Fig. \ref{fig3+}. First, we prove this lemma when $\alpha>\theta$. According to the Fig. \ref{fig3+}, 
$$|DC|=|DB|=r,\hspace{1.5cm}\beta:=\angle CDB=\alpha-\theta,\hspace{1.5cm}|AB|=\frac{1}{\sin\alpha}.$$
\par
Moreover, the distance between the point $C$ and the boundary of the corridor is equal to $r-r\cos\beta$. This implies that,
$$\sin(\angle CAB)>(r-r\cos\beta)\sin\alpha=(r-r\cos(\alpha-\theta))\sin\alpha.$$
If we denote by $L$ the length of the longest segment in this corridor with one endpoints at $C$, then:
$$L<\frac{d}{\sin(\angle CAB)}+2|AB|<\frac{d}{(r-r\cos(\alpha-\theta))\sin\alpha}+\frac{2}{\sin\alpha}.$$
To complete the proof for this case, let $L_0=\frac{d}{(r-r\cos(\alpha-\theta))\sin\alpha}+\frac{2}{\sin\alpha}$.\par
In case $\alpha<\theta$, the corridor is tangent to the dispersing components on left and right sides of $\partial S'$. The proof in this case is very similar to the previous one, and it results at:
$$L_0=\frac{d}{(r-r\cos(\theta-\alpha))\cos\alpha}+\frac{2}{\cos\alpha},$$
when $\alpha<\theta$.\qed\par
\begin{figure}[h]
	\centering
	\includegraphics[width=9cm]{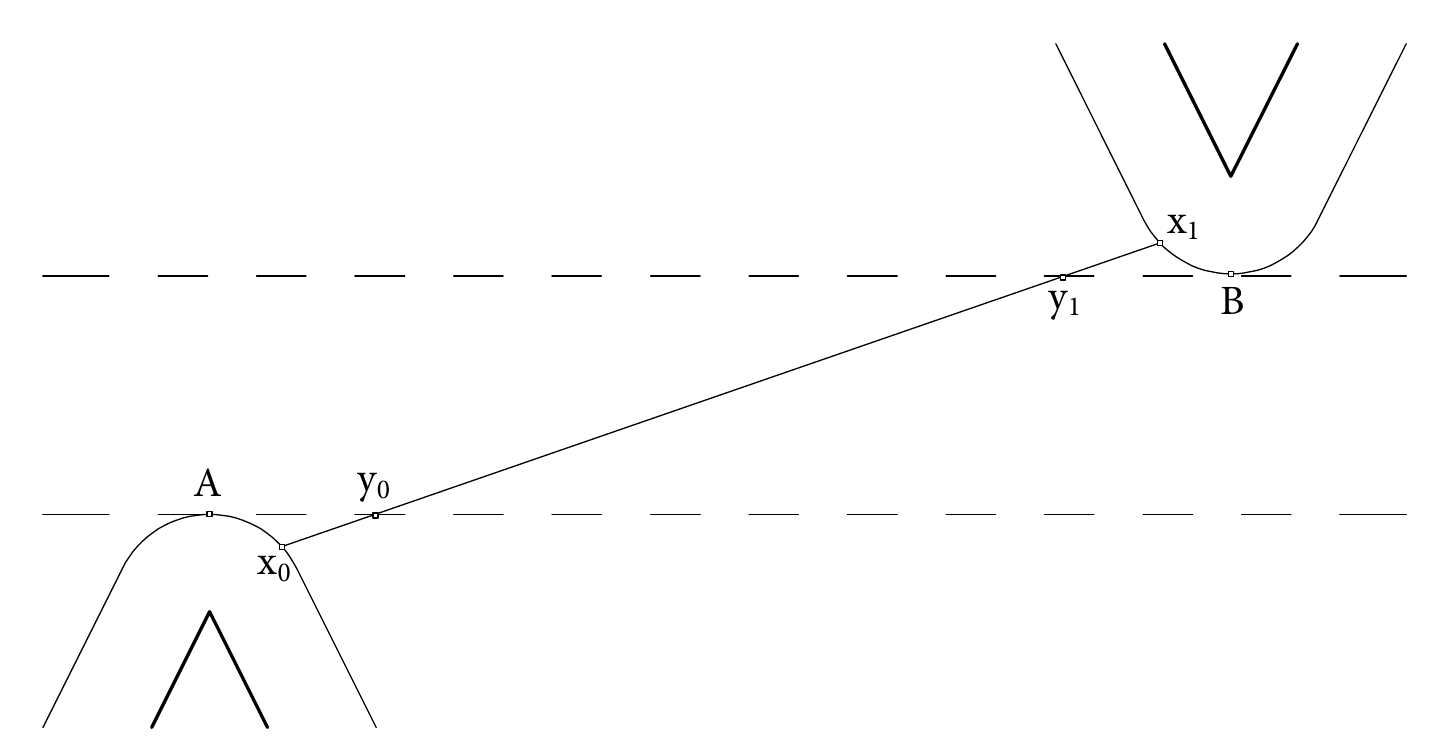}\caption{Intersection of a finite segment $[x_0,x_1]$ with the boundary of $C_h$.}\label{lem4}
\end{figure}

Let $L_h$ and $L_v$ denote the parameter $L_0$ from Lemma \ref{lem3}, respectively, for the horizontal and vertical corridors $C_h$ and $C_v$. Then,
$$L_h=\frac{d_h}{r-r\sin\theta}+2,\hspace{1cm}L_v=\frac{d_v}{r-r\cos\theta}+2.$$
Consider the segment $[x_0,x_1]$ is in $C_h$ or $C_v$ such that $|x_0-x_1|>\max\{L_h,L_v\}$. If $y_0$ and $y_1$ are the intersection points of this segment with the boundary of the corridor (Fig. \ref{lem4}), then
\begin{equation}\label{lem4-1}
	|x_0-y_0|<1,\hspace{1in}|x_1-y_1|<1.
\end{equation}
From Lemma \ref{lem3}, we know that $x_0$ and $x_1$ are on dispersing components of $\partial S'$. Let $A$ and $B$ denote the tangent points to the boundaries of the corridor of those dispersing components which contain the points $x_0$ and $x_1$, respectively. (Fig. \ref{lem4})
\begin{Lemma}\label{lem4+}
	If $[x_0,x_1]$ is a segment in $C_h$ or $C_v$ such that $|x_0-x_1|>\max\{L_h,L_v\}$, then there exists an $\epsilon>0$ such that:
	$$|x_0-A|<\epsilon,\hspace{1in}|x_1-B|<\epsilon.$$
	Moreover, $|x_0-A|$ and $|x_1-B|$ tend to zero when $|x_0-x_1|\rightarrow\infty$.
\end{Lemma}
\proof
It is easy to see that $\epsilon=r(\frac{\pi}{2}-\theta)$. Moreover, the calculation in the proof of Lemma \ref{lem3} suggests that $|x_0-A|$ and $|x_1-B|$ tend to zero when $|x_0-x_1|\rightarrow\infty$.\qed\par
The inequalities in (\ref{lem4-1}) and lemmas \ref{lem3} and \ref{lem4+} guarantee that asymptotics of the probability of free motion vectors in horizontal and vertical corridors of the periodic Ehrenfests' Wind-Tree model is the same as asymptotics of the probability of free motion vectors in corridors of the periodic Lorentz gas \cite{Bleher92}. According to the result of proposition 4.2 in \cite{Bleher92}, the expressions for the corresponding probabilities in the limit $L\rightarrow\infty$ are:
\begin{equation}\label{asymH}
	\begin{split}
		\mathbf{Pr}_h(L):= & Probability\ of\ \{r(X)\ in\ C_h\ such\ that\ |r(X)|>L\} \\ = & \frac{Z^{-1}d_h^2}{L^2}+O(L^{-5/2}),
	\end{split}
\end{equation}
and
\begin{equation}\label{asymV}
	\begin{split}
		\mathbf{Pr}_v(L):= & Probability\ of\ \{r(X)\ in\ C_v\ such\ that\ |r(X)|>L\} \\ = & \frac{Z^{-1}d_v^2}{L^2}+O(L^{-5/2}).
	\end{split}
\end{equation}
\section{Asymptotics of $\nu$ in Corridors of Type II}\label{sec4}
Our main result is concerned with oblique corridors of type II, where the boundary of the corridor contains flat components of scatterers. We will show that the existence of such corridors results in stronger superdiffusive regimes than the one in the Lorentz gas.\par
\begin{figure}[h]
	\centering
	\includegraphics[width=7cm]{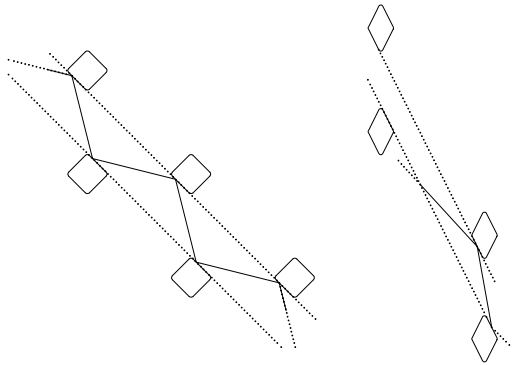}\caption{Trajectories with reflections on neutral parts in a corridor of type II when $\theta=\pi/4$ (left) and $\theta\neq\pi/4$ (right).}\label{figNOI}
\end{figure}

From Lemma \ref{lem3}, we can expect oblique corridors of type I to have the same diffusive properties as the horizontal and vertical corridors. To reduce the volume of calculations, we will consider a physical periodic Ehrenfests' Wind-Tree model without oblique corridors of type I. Moreover, if there exists an oblique corridor of type II where $\theta\neq\frac{\pi}{4}$, then the scatterers are distributed over the boundary of this corridor in a way that any trajectory in this corridor is likely to leave the corridor after a few reflections on neutral parts of scatterers (See Fig. \ref{figNOI}).\par
By Lemma \ref{lem2}, one can calculate the minimum value of $a$ such that there is no oblique corridor of type I when $\theta=\frac{\pi}{4}$. Since $\tan(\theta)=1$, we only need to consider cases 2 and 4 in Lemma \ref{lem2}. Moreover, under this assumption, the upper bounds of $a$ in these two cases will be the same and equal to:
$$\frac{n+m-m\lceil\frac{n}{m}\rceil}{\sqrt{2}n},$$
where $0<m<n$ are integers. Therefore, if $\theta=\frac{\pi}{4}$ and
$$a\geq Max_{\{(n,m)|\ n>m,\ n,m\in\mathbf{N}\}}\left\{\frac{n+m-m\lceil\frac{n}{m}\rceil}{\sqrt{2}n}\right\}=\frac{\sqrt{2}}{4},$$
then there are no oblique corridors of type I.\par
In the rest of this paper, we consider physical periodic Ehrenfests' Wind-Tree model with a moving disk of radius $0<r<\frac{\sqrt{2}}{8}$ where $\theta=\frac{\pi}{4}$ and $\frac{\sqrt{2}}{4}\leq a<\frac{\sqrt{2}}{2}-2r$.\par
Thus, our model has two corridors of type I ($C_h$ and $C_v$) and two corridors of type II with directions parallel to $y=x$ and $y=-x$. Denote these oblique corridors of type II by $C_o^+$ and $C_o^-$, respectively.\par
\begin{figure}[h]
	\centering
	\includegraphics[width=11cm]{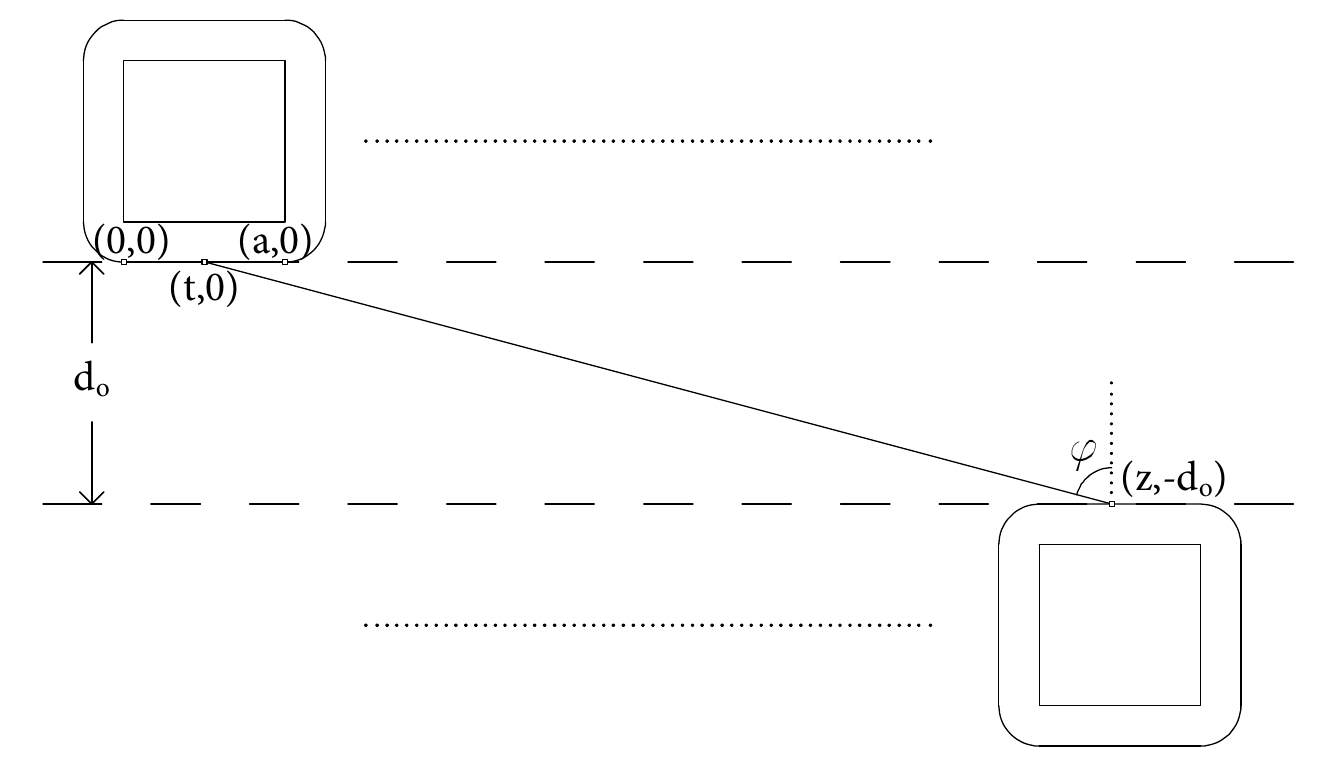}\caption{A long collision free segment in an oblique corridor of type II.}\label{fig32}
\end{figure}

To estimate the tail of the distribution $\nu$ in this model along corridors $C_o^+$ and $C_o^-$, we need to calculate the asymptotics of
\begin{equation*}
	\begin{split}
		\mathbf{Pr}_o(L):= & Probability\ of\ \{r(X)\ in\ C_o^+\ such\ that\ |r(X)|>L\} \\ = & Probability\ of\ \{r(X)\ in\ C_o^-\ such\ that\ |r(X)|>L\}
	\end{split}
\end{equation*}
when $L\rightarrow\infty$.
\begin{proposition}\label{prop1}
	When $L\rightarrow\infty$,
	\begin{equation}\label{asymO}
		\mathbf{Pr}_o(L)=aZ^{-1}\left(\frac{d_o}{L}\right)^2+O(L^{-3}),
	\end{equation}
	where $d_o$ is the width of oblique corridors $C_o^+$ or $C_o^-$.
\end{proposition}
\proof
According to the Fig. \ref{fig32}, in the oblique corridors of type II of the physical periodic Wind-tree model: 
$$w:=\tan(\varphi)=\frac{z-t}{d_o}\ \Longrightarrow\ \cos(\varphi)d\varphi=\frac{dz}{d_o(1+w^2)^{3/2}}.$$
Therefore, the probability density function $f(.)$ of the distribution of $z$ with respect to the Liouville measure $Z^{-1}\cos(\varphi)d\varphi ds$ along these corridors is equal to:
$$f(z)=\int_0^a\frac{Z^{-1}}{d_o(1+w^2)^{3/2}}dt=Z^{-1}d_o^2\int_0^a\frac{dt}{(d_o^2+(z-t)^2)^{3/2}}$$
$$\simeq Z^{-1}d_o^2\int_0^a\frac{dt}{(z-t)^3}=\frac{Z^{-1}d_o^2}{2}\left(\frac{1}{(z-a)^2}-\frac{1}{z^2}\right)=\frac{aZ^{-1}d_o^2}{z^3}(1+O(z^{-1})).$$
Hence, the probability of $|r(X_0)|>L$ is given by:
\begin{equation}\label{coe.2}
	\begin{split}
		\mathbf{Pr}_o(L)=& 2\int_L^\infty f(z)dz
		=2\int_L^\infty\frac{aZ^{-1}d_o^2}{z^3}(1+O(z^{-1}))dz\\
		=& aZ^{-1}\left(\frac{d_o}{L}\right)^2+O(L^{-3}).
	\end{split}
\end{equation}\qed\par
The factor $2$ in the expression for $\mathbf{Pr}_o(L)$ in the equation (\ref{coe.2}) appears because there are two opposite directions in each corridor.
\begin{proposition}\label{prop2}
	$\langle r(X)\rangle=0$, $\langle|r(X)|\rangle<\infty$, and $\langle|r(X)|^2\rangle=\infty$. Moreover, if $\phi_R(x)=|x|^2$ when $|x|<R$ and it is zero otherwise, then
	$$\langle\phi_R(r(X))\rangle=Const.\ln(R)+O(1),$$
	when $R\rightarrow\infty$.
\end{proposition}
\proof
We have:
$$\langle r(X)\rangle=\int r(X)\mu(d\varphi ds)=0,$$
because $\nu$ is symmetric (Proposition \ref{nu-sym}). If $L$ is sufficiently large, then it follows from (\ref{asymH}-\ref{asymO}) that:
$$\langle|r(X)|\rangle=\int|r(X)|\mu(d\varphi ds)=\int_{|r|\leq L}|r(X)|\mu(d\varphi ds)+\int_{|r|>L}|r(X)|\mu(d\varphi ds)$$
$$=\int_{|r|\leq L}|r(X)|\mu(d\varphi ds)+\int_{\substack{|r|>L\\in\ C_o^+}}|r(X)|\mu(d\varphi ds)+\int_{\substack{|r|>L\\in\ C_o^-}}|r(X)|\mu(d\varphi ds)$$
$$+\int_{\substack{|r|>L\\in\ C_h}}|r(X)|\mu(d\varphi ds)+\int_{\substack{|r|>L\\in\ C_v}}|r(X)|\mu(d\varphi ds)$$
$$<L+2\int_L^\infty\frac{2aZ^{-1}d_o^2}{x^2}dx+\int_L^\infty\frac{2Z^{-1}d_h^2}{x^2}dx+\int_L^\infty\frac{2Z^{-1}d_v^2}{x^2}dx<\infty.$$
Finally, for $R\rightarrow\infty$:
\begin{equation*}
	\begin{split}
		\langle\phi_R(r(X))\rangle=&\int_{|r|\leq L<R}|r(X)|^2\mu(d\varphi ds)+\int_{L<|r|<R}|r(X)|^2\mu(d\varphi ds)\\
		=& M_0+\int_{\substack{L<|r|<R\\in\ C_o^+}}|r(X)|^2\mu(d\varphi ds)+\int_{\substack{L<|r|<R\\in\ C_o^-}}|r(X)|^2\mu(d\varphi ds)\\
		&\ \ \ \ +\int_{\substack{L<|r|<R\\in\ C_h}}|r(X)|^2\mu(d\varphi ds)+\int_{\substack{L<|r|<R\\in\ C_v}}|r(X)|^2\mu(d\varphi ds)\\
		=& M_0+2\int_L^R\frac{2aZ^{-1}d_o^2}{x}dx+\int_L^R\frac{2Z^{-1}d_h^2}{x}dx+\int_L^R\frac{2Z^{-1}d_v^2}{x}dx\\
		=& Const.\ln(R)+O(1).
	\end{split}
\end{equation*}
Therefore, $\langle|r(X)|^2\rangle=\infty$.\qed
\section{Statistical Properties of Ehrenfests' Wind-Tree Models}
In the previous section, we studied the asymptotic behavior of the distribution $\nu$ of the free motion vector $r(X)=x(TX)-x(X)$. Now, we will consider the joint distribution $\nu_n$ of vectors $r(X)$ and $r(T^nX)$ with respect to the Liouville measure $\mu$ in order to estimate,
$$\langle|r(X)||r(T^nX)|\rangle,$$
for any $n\neq0$.
\begin{proposition}\label{prop5}
	The distribution $\nu_n$ is invariant with respect to the transformation $(r_1,r_2)\mapsto-(r_2,r_1)$.
\end{proposition}
The proof is identical to the proof of Proposition 5.1 in \cite{Bleher92}. Moreover, the Proposition \ref{prop5} implies,
$$(|r(X_0)|,|r(T^nX_0)|)\stackrel{d}{=}(|r(T^nX_0)|,|r(X_0)|).$$
\subsection{Estimation of Correlations of Free Motion Vectors}
Let $r(X_0)=[x_0,x_1]$ be in $C_o^+$ or $C_o^-$ and $x_1$ be on a neutral component of $\partial S'$ belongs to the boundary of that corridor. Then quantities $|r(X_0)|$ and $|r(TX_0)|$ satisfy the relation:
\begin{equation}\label{relation o}
	|r(TX_0)|=|r(X_0)|+O(1).
\end{equation}
More generally, when endpoints of segments  $r(T^{i-1}X_0)=r(X_{i-1})=[x_{i-1},x_i]$ for $i=1,\dots,n$ belong to neutral components of $\partial S'$ which are in the boundary of a corridor of type II, then:
\begin{equation}\label{relation o G}
	|r(X_i)|=|r(T^iX_0)|=|r(X_0)|,
\end{equation}
for $i=1,2,\dots,n-1$, and
\begin{equation}\label{relation o G+}
	|r(X_0)|\leq|r(X_n)|=|r(T^nX_0)|<|r(X_0)|+1.
\end{equation}
Denote the expected value of a function along the corridors $C_o^+$ or $C_o^-$ by $\langle.\rangle_o$.
\begin{proposition}\label{prop6}
	Let all $x_i$ of segments  $r(T^iX)=[x_i,x_{i+1}]$ for $i=0,1,\dots,n$ belong to neutral components of $\partial S'$ in the boundary of $C_o^+$ or $C_o^-$, then
	$$\langle|r(X)||r(T^{n-1}X)|\rangle_o=\frac{\langle|r(X)|^2\rangle_o}{n}.$$
\end{proposition}
\proof
First, we find the probability of trajectories with $n$ consecutive reflections on neutral parts of $\partial S'$ in the boundary of a corridor of type II.\par
\begin{figure}[h]
	\centering
	\includegraphics[width=11cm]{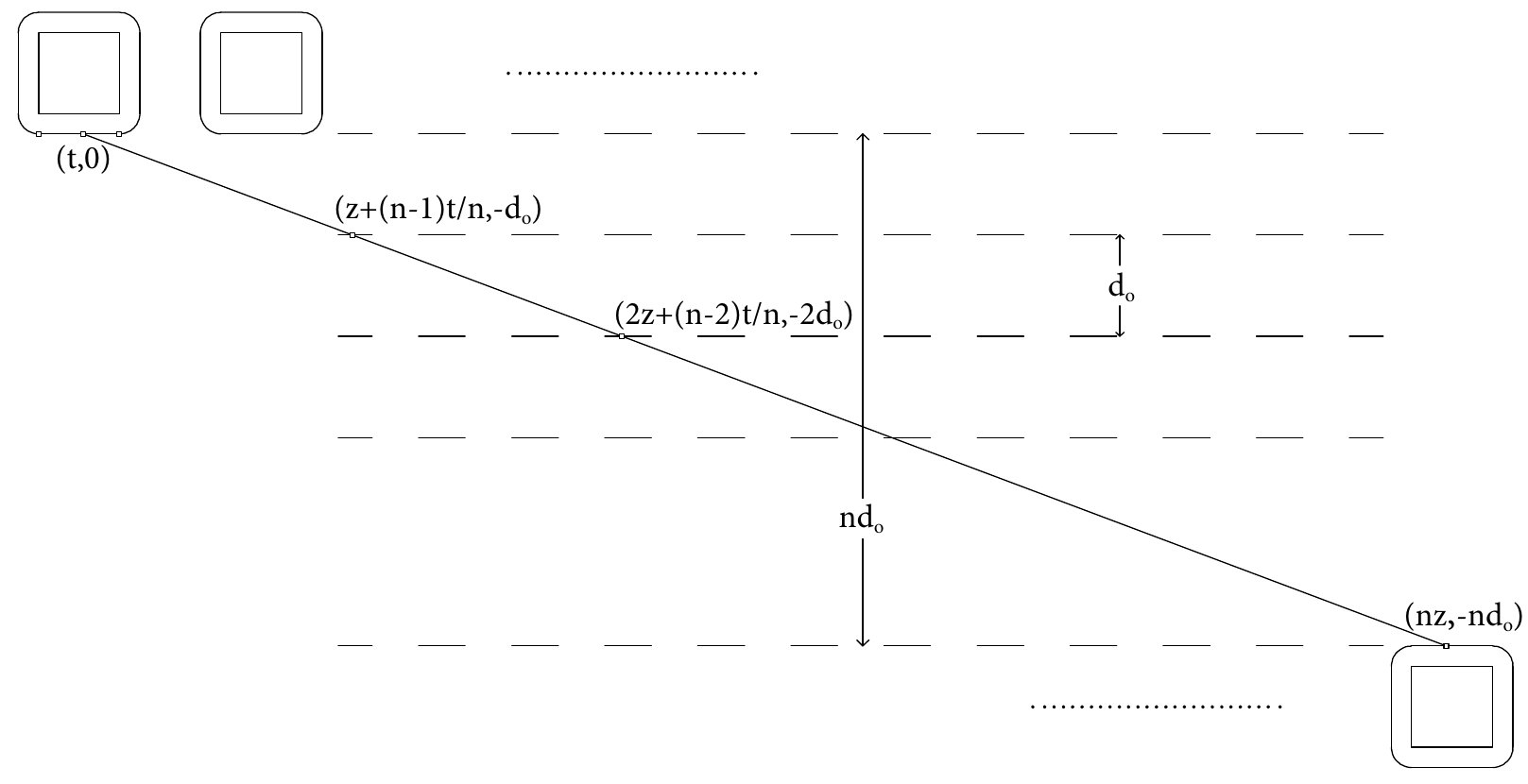}\caption{Segments with reflections on neutral components of $\partial S'$ in boundaries of a corridor of type II.}\label{fig33}
\end{figure}
It follows from Fig. \ref{fig33} that,
$$w=\tan(\varphi)=\frac{nz-t}{nd_o}\ \Longrightarrow\ \cos(\varphi)d\varphi=\frac{dz}{nd_o(1+w^2)^{3/2}}.$$
Then,
\begin{equation}\label{f_n}
	\begin{split}
		f_n(z):=&\int_0^a\frac{Z^{-1}}{nd_o(1+w^2)^{3/2}}dt=Z^{-1}n^2d_o^2\int_0^a\frac{dt}{(n^2d_o^2+(nz-t)^2)^{3/2}}\\
		\simeq& Z^{-1}n^2d_o^2\int_0^a\frac{dt}{(nz-t)^3}=\frac{Z^{-1}n^2d_o^2}{2}\left(\frac{1}{(nz-a)^2}-\frac{1}{n^2z^2}\right)\\
		=&\frac{aZ^{-1}d_o^2}{nz^3}(1+O(z^{-1})).
	\end{split}
\end{equation}
The equation (\ref{f_n}) and Proposition \ref{prop1} show that,
\begin{equation}\label{f_n1}
	f_n(z)=\frac{1}{n}f(z).
\end{equation}
Let $P_n(L)$ denote the probability of orbits with segments $[x_i,x_{i+1}]$ for $i=0,1,\dots,n$ such that all $x_i$ are on the neutral parts in the boundary of a corridor of type II and $|x_{i+1}-x_i|>L$. Then, it follows from (\ref{f_n1}) and the proof of Proposition \ref{prop1} that,
$$P_n(L)=\frac{1}{n}\mathbf{Pr}_o(L).$$
Moreover, the equation (\ref{relation o G}) implies,
$$|r(X)||r(T^{n-1}X)|=|r(X)|^2.$$
Therefore,
$$\langle|r(X)||r(T^{n-1}X)|\rangle_o=\frac{\langle|r(X)|^2\rangle_o}{n}.$$\qed\par
We will now estimate $\langle|r(X)||r(T^{n-1}X)|\rangle_o$ where some $x_i$ for $i=0,1,2,\dots,n$ belongs to the dispersing parts of $\partial S'$. Let consider the case that there is a reflection off a dispersing component and it is followed by $n$ consecutive reflections on neutral parts in the boundary of a corridor of type II. The proof of Proposition \ref{prop6}, particularly the equation (\ref{f_n}), shows that:
$$\langle|r(X)||r(T^nX)|\rangle_o=\frac{\langle|r(X)||r(TX)|\rangle_o}{n},$$
where $x_1$ in $r(X)=[x_0,x_1]$ is on a dispersing component. Then, proposition 5.3 in \cite{Bleher92} implies that $\langle|r(X)||r(TX)|\rangle_o=Const$. Therefore, in this case
\begin{equation}\label{discor}
	\langle|r(X)||r(T^nX)|\rangle_o=\frac{\langle|r(X)||r(TX)|\rangle_o}{n}=\frac{O(1)}{n}.
\end{equation}
More generally, when an orbit has only reflections on neutral parts in the boundary of a corridor of type II, the angle between its segments and the length of its segments remain the same, while reflections on dispersing components change these angles and they also change lengths of free paths. The average of these changes (correlations) is estimated in theorem 5.6 in \cite{Bleher92} when an orbit has only reflections on dispersing components. The results of \cite{Bleher92} and Proposition \ref{prop6} show that when an orbit has more than one reflection on dispersing components, the correlations is smaller than the correlations in the case where there is only one reflection on a dispersing part when the particle goes along a corridor of type II. Hence, from the equation (\ref{discor}) and Proposition \ref{prop6}, the correlation function in a corridor of type II satisfies the following relation:
\begin{equation}\label{corO}
	\langle(r(X),r(T^nX))\rangle_o\simeq\langle|r(X)||r(T^nX)|\rangle_o=\frac{\langle|r(X)|^2\rangle_o+O(1)}{n+1},
\end{equation}
since $r(X)$ and $r(T^nX)$ have almost the same directions in that corridor of type II (i.e. $(r(X),r(T^nX))\simeq|r(X)||r(T^nX)|$).\par
To better understand this result observe that two neighboring boundary components of a straight segment in the boundary of any scatterer are arcs of a circle with the same radius. Therefore their combined influence on correlations is the same as in the Lorentz gas \cite{Bleher92}. Indeed, by putting them together we get circular part of scatterer which influences the passage of a particle through the corridor. The fact that dispersing components are apart on the length of a neutral component clearly contributes only to a constant factor in the estimate of correlations.
\subsection{Statistical Behavior of Trajectories in Discrete Dynamics}
From the result of previous sections, one can expect that the physical periodic Ehrenfests' Wind-Tree models will have more regimes of diffusion than the Lorentz gas. To show that we need to estimate $\langle|x(T^nX)-x(X)|^2\rangle$ when $n\rightarrow\infty$. This value can be written as:
\begin{equation}\label{1steq}
	\begin{split}
		\langle|x(T^nX)-x(X)|^2\rangle&=\sum_{i=0}^{n-1}\sum_{j=0}^{n-1}\langle(r(T^iX),r(T^jX))\rangle\\
		&=n\langle|r(X)|^2\rangle+2\sum_{j=1}^{n-1}(n-j)\langle(r(X),r(T^jX))\rangle.
	\end{split}
\end{equation}
It follows from (\ref{corO}) that,
\begin{equation}\label{suminf}
	\begin{split}
		\lim_{n\rightarrow\infty}\sum_{k=1}^n\langle(r(X),r(T^kX))\rangle_o
		&\simeq\lim_{n\rightarrow\infty}\sum_{k=1}^n\frac{\langle|r(X)|^2\rangle_o+O(1)}{k+1}\\ & =\lim_{n\rightarrow\infty}[(\ln(n)-1)\langle|r(X)|^2\rangle_o+O(\ln n)].
	\end{split}
\end{equation}
Moreover, Proposition \ref{prop2} implies:
\begin{equation}\label{<>O}
	\begin{split}
		\langle\phi_R(r(X))\rangle_o & =\int_{\substack{|r|\leq L<R\\in\ C_o^+}}|r(X)|^2\mu(d\varphi ds)+\int_{\substack{L<|r|<R\\in\ C_o^+}}|r(X)|^2\mu(d\varphi ds)\\
		& =m_0+2\int_L^R\frac{aZ^{-1}d_o^2}{x}dx=\frac{2ad_o^2}{Z}\ln R+O(1),
	\end{split}
\end{equation}
when $R\rightarrow\infty$. It follows from equations (\ref{suminf}) and (\ref{<>O}) that,
\begin{equation}\label{lnn2}
	\sum_{k=1}^N\langle(r(X),r(T^kX))\rangle_o=\frac{2ad_o^2}{Z}(\ln N)^2+O(\ln N).
\end{equation}
Therefore the main contribution to correlations is made by orbits propagating in corridors of type II and reflecting only off neutral components of the boundary.\par
Consider a physical periodic Wind-Tree model. Then for any positive radius $r>0$ of the moving particle we conjecture that the following statement is correct. Already existing methods \cite{Bleher92,Bun81,BSC91} are more than enough for its proof, which does not require any new ideas besides those presented in our paper.
\begin{Conjecture}
	Let the distribution of $X\in\Lambda$ be the Liouville measure $\mu(d\varphi ds)=Z^{-1}\cos(\varphi)d\varphi ds$. If
	\begin{equation}\label{main}
		\xi=\lim_{n\rightarrow\infty}\frac{x(T^nX)-x(X)}{g(n)},
	\end{equation}
	then $\xi=(\xi_1,\xi_2)$ is a Gaussian random variable with zero mean where,
	\begin{enumerate}
		\item $g(n)=\sqrt{n}$ in case of finite horizon.
		\item $g(n)=\sqrt{n\ln n}$ in case of infinite horizon but without corridors of type II.
		\item $g(n)=\sqrt{n}\ln n$ in case of infinite horizon with presence of type II corridors.
	\end{enumerate}
	Moreover, in item 3, when we consider the physical periodic Wind-Tree model presented in Section \ref{sec4}, the covariance matrix is given by
	$$\left[\begin{array}{l}
	D_{11}\hspace{0.6cm}D_{12}\\
	D_{21}\hspace{0.6cm}D_{22}
	\end{array}\right]=\left[\begin{array}{l}
	\frac{2ad_o^2}{|\partial S'|}\hspace{0.7cm}0\\
	\ \ 0\hspace{0.8cm}\frac{2ad_o^2}{|\partial S'|}
	\end{array}\right].$$
\end{Conjecture}
An outline of the proof in case of infinite horizon with presence of corridors of type II is as follows. According to Proposition \ref{nu-sym}, the probability distribution $\nu$ is symmetric, i.e. $\langle\xi\rangle=0$. To find the covariance matrix, we need to calculate its components in corridors. There are four corridors in the physical periodic Ehrenfests' Wind-Tree model: $C_h$, $C_v$, $C_o^+$, and $C_o^-$. It follows from Section \ref{sec3} that asymptotics of the probability distribution in $C_h$ and $C_v$ are the same as those in corridors of the Lorentz gas. Thus, with normalization $\sqrt{n}\ln(n)$ used in (\ref{main}) and the result from \cite{Bleher92}, we will have $\langle\xi_i\xi_j\rangle=0$ for $i,j=1,2$ when $\xi$ is in $C_h$ or $C_v$. This means, $D_{ij}=\langle\xi_i\xi_j\rangle$ only depends on the value $\langle\xi_i\xi_j\rangle_o$ where $\xi$ is in $C_o^+$ or $C_o^-$. Let $\xi$ be in $C_o^+$, then:
\begin{equation}\label{++}
	\xi_1\simeq\xi_2,
\end{equation}
and,
\begin{equation}\label{+}
	\langle\xi_1\xi_2\rangle_o=\langle\xi_1^2\rangle_o.
\end{equation}
Similarly, when $\xi$ is in $C_o^-$:
\begin{equation}\label{--}
	\xi_1\simeq-\xi_2,
\end{equation}
and,
\begin{equation}\label{-}
	\langle\xi_1\xi_2\rangle_o=-\langle\xi_1^2\rangle_o.
\end{equation}
From (\ref{+}) and (\ref{-}), it follows that:
$$D_{12}=D_{21}=sum\ of\ values\ of\ \langle\xi_1\xi_2\rangle_o\ along\ C_o^+\ and\ C_o^-=0$$
On the other hand, (\ref{++}) and (\ref{--}) imply that,
\begin{equation}\label{d11}
	\xi_i^2\simeq\frac{|\xi|^2}{2},
\end{equation}
for $i=1,2$.\par
By making use of (\ref{1steq}), (\ref{lnn2}), and (\ref{d11}) and Proposition \ref{prop2}, one gets:
\begin{equation*}
	\begin{split}
		D_{11}=D_{22}=&sum\ of\ values\ of\ \langle\xi_1^2\rangle_o\ along\ C_o^+\ and\ C_o^-=2\langle\xi_1^2\rangle_o=\langle|\xi|^2\rangle_o\\= & \langle\lim_{n\rightarrow\infty}\frac{|x(T^nX)-x(X)|^2}{n(\ln(n))^2}\rangle_o=\lim_{n\rightarrow\infty}\frac{\langle|x(T^nX)-x(X)|^2\rangle_o}{n(\ln(n))^2}\\= & \lim_{n\rightarrow\infty}\frac{n\langle|r(X)|^2\rangle_o+2\sum_{j=1}^{n-1}(n-j)\langle(r(X),r(T^jX))\rangle_o}{n(\ln(n))^2}\\= & \lim_{n\rightarrow\infty}\frac{2\sum_{j=1}^{n-1}\langle(r(X),r(T^jX))\rangle_o}{(\ln(n))^2}=\frac{4ad_o^2}{Z}=\frac{2ad_o^2}{|\partial S'|}.
	\end{split}
\end{equation*}
\subsection{Statistical Behavior of Trajectories in Continuous-Time Dynamics}
In this section, we will present an analogous formula to (\ref{main}) for continuous-time dynamics. Let $t_n$ be the time of the $n^{th}$ reflection of the trajectory $x(t)$. Then, $x(t_n)=x_n=x(X_n)$. Since the path of the particle between reflections is a line segment and the particle moves with velocity $1$ along it,
$$t_{n+1}-t_n=|x_{n+1}-x_n|=|r(X_n)|,$$
and,
$$t_n=\sum_{i=0}^{n-1}|r(T^iX_0)|=\sum_{i=0}^{n-1}|r(X_i)|.$$
For an arbitrary initial condition $X_0$, the ergodic theorem guarantees that,
$$\lim_{n\rightarrow\infty}\frac{t_n}{n}=\lim_{n\rightarrow\infty}\frac{\sum_{i=0}^{n-1}|r(X_i)|}{n}=\langle|r(X_0)|\rangle.$$
From Proposition \ref{prop2}, $\eta:=\langle|r(X_0)|\rangle<\infty$. Therefore,
$$\lim_{n\rightarrow\infty}\frac{t_n(\ln t_n)^2}{n(\ln n)^2}=\eta,$$
and
$$\lim_{n\rightarrow\infty}\frac{t_{n+1}(\ln t_{n+1})^2}{t_n(\ln t_n)^2}=1.$$
Let $t_n\leq t<t_{n+1}$. It follows from lemma 7.1 in \cite{Bleher92} and last two equations that,
\begin{equation}\label{continuous}
	\lim_{n\rightarrow\infty}\frac{x_t-x_0}{\sqrt{t}\ln t}=\lim_{n\rightarrow\infty}\frac{x_t-x_0}{\sqrt{t_n}\ln t_n}=\eta^{-1/2}\lim_{n\rightarrow\infty}\frac{x_t-x_0}{\sqrt{n}\ln n}.
\end{equation}
The assumption $t_n\leq t<t_{n+1}$ implies that $x_t$ belongs to the segment $[x_n,x_{n+1}]$. Hence, $|x_t-x_n|\leq|r(X_n)|$. Thus almost surely,
$$\lim_{n\rightarrow\infty}\frac{x_t-x_n}{\sqrt{n}\ln n}=0.$$
Therefore,
$$\lim_{n\rightarrow\infty}\frac{x_t-x_0}{\sqrt{n}\ln n}=\lim_{n\rightarrow\infty}\frac{x_n-x_0}{\sqrt{n}\ln n}=\xi.$$
where $\xi$ is the same as in Conjecture. Then, it follows from (\ref{continuous}) that,
$$\lim_{t\rightarrow\infty}\frac{x_t-x_0}{\sqrt{t}\ln t}=\frac{\xi}{\sqrt{\eta}}.$$
This expression is analogous to (\ref{main}) when we consider the continuous-time dynamics, and it shows that there is a new superdiffusive regime in the physical periodic Ehrenfests' Wind-Tree model where the diffusion coefficient $D(t)\sim(\ln t)^2$.
\section*{Acknowledgments}
We thank R. Feres for useful suggestions and for pointing to the paper \cite{CFZ18}. This work was partially supported by the NSF grant DMS-1600568.

\end{document}